# ON AGEING PROPERTIES OF LIFETIME DISTRIBUTIONS

Anakha K K & V M Chacko

*Department of Statistics, St. Thomas' College (Autonomous),*

*Thrissur, Kerala-680 001, India*

Email: anakhaappu@gmail.com, chackovm@gamil.com

**Abstract**

A reasonable segment of reliability theory is perpetrated to the study of failure rates, their properties, connections and applications. The present study focused on failure rate distributions and their shape properties. Failure rates of various generalizations of Lindley distribution are reviewed and the distinguishable role of DUS transformation is also discussed.

1. **Introduction**

Lifetime distributions offer a clear perception of reliability theory and survival analysis as the characteristic terms such as failure rate function and mean residual life function gives comprehensive accounts through real-time data modelling. Over the preceding decades, a considerable thorough study has been done on lifetime distributions. The failure rate has several applications in diverse areas and it is known by distinctive names. In actuaries, failure rates are used to compute mortality tables so-called "force of mortality" (Steffensen(1930)). In statistics, the reciprocal of failure rate for the Normal distribution is known "Mills ratio". While in extreme value theory, failure rate plays an important role in the so-called "intensity function" (Gumbel(1958)), which is also called the "hazard rate" in reliability theory.

The term'ageing' plays a crucial role in reliability where we generally specify the lifetime distribution through the ageing properties. According to Lai *et.al* (2006), the term 'no ageing' refers to the age of the component not affect the distribution of residual life length. In a probabilistic contest, 'positive ageing' ('averse ageing') describes the circumstance in which residual lifetime decreases as the age of the component increases. Negative ageing ('beneficial ageing'), on the other hand, has the opposite effect on residual lifetime.

These characteristic terms have been used to discuss ageing qualities such as increasing failure rate (IFR), decreasing failure rate (DFR), and so on. Other categories include 'increasing failure rate on average' (IFRA), 'decreasing failure rate on average' (DFRA), 'new better than used' (NBU), 'new better than used in expectation (NBUE), and 'decreasing mean residual life' (DMRL). For a complete account see Barlow and Proschan (1975) and Hollander and Proschan (1984). Table 1 gives the abbreviation of some ageing classes that is used to examine in this paper.

**Table 1**

| Abbreviation | Name of class |
|---|---|
| IFR | Increasing failure rate |
| IFRA | Increasing failure rate average |
| NBU | New better than used |
| NBUE | New better than used in expectation |
| IMRL | Increasing mean residual life |
| HNBUE | Harmonic new better than used |
| DFR | Decreasing failure rate |
| DFRA | Decreasing failure rate average |
| NWU | New worse than used |
| NWUE | New worse than used in expectation |
| DMRL | Decreasing mean residual life |
| $L$ | Class with ageing property based on Laplace transform |
| BFR | Bathtub-shaped failure rate |
| UFR | Upside-down bathtub shaped failure rate |
| NBU-$t_0$ | New better than used of age $t_0$ |
| NWU-$t_0$ | New worse than used of age $t_0$ |

This paper mainly focused on the basic ageing properties of lifetime distributions. The rest of the paper is organized as follows. Definitions of failure rates and mean residual life and some of the properties of ageing classes have been discussed in section 2. An outline of IFR, DFR and their properties has been included in section 3. A summary of roller-coaster shaped failure rates is given in section 4. In section 5, various generalizations of Lindley distribution and their failure rates have been discussed. A study on DUS transformed distributions and their properties are contributed in section 6. Additionaly, DUS transformation on mixture of Exponential-Weibull distribution is proposed. Finally, a conclusion has been given in section 7.

## 2. Failure rates and Mean residual life

The concept of failure rate has been beneficial in instructing an imperfect repair model when proper maintenance is performed. The probability density function (*pdf*) and cumulative distribution function (*cdf*) of the lifespan random variable *T*, respectively, are denoted by *f(t)* and *F(t)*. The failure rate function is $h(t) = \dfrac{f(t)}{\overline{F}(t)}$ where $\overline{F}(t) = 1 - F(t)$, t $\geq$ 0. In lifetime experiments, the average additional lifetime of components when they remain up to time t is called mean residual life(MRL) function (Gupta and Kirmani (2000)). The MRL is characterised by $M(t) = E(T - t | T > t), t \geq 0$.

Since the MRL is closely connected to the failure rate function, some studies examined the monotonicity of MRL concerning the monotonicity of the failure rate function. For a detailed account of MRL, it is recommended to read Guess and Proschan (1988).

The basic ageing properties are related by the following illustrations (Deshpande *et al*. (1986) and Rao (1992))

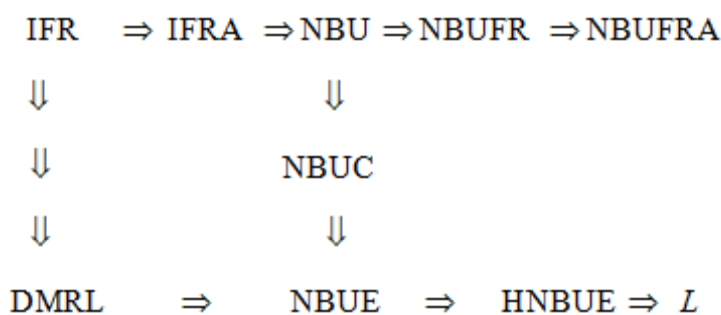

Unnikrishnan and Vineshkumar (2011) have studied ageing properties in terms of quantile function.

**Definition 2.1** (Lai *et al.* (2006))

Assume that the failure rate *h(t)* is a real-valued differentiable function $h(t): R^+ \to R^+$. Then the cumulative distribution function *F(t)* is said to be

1. Increasing failure rate (IFR) if $h'(t) > 0$ for all t;
2. Decreasing failure rate (DFR) if $h'(t) < 0$ for all t;
3. Bathtub shaped failure rate (BFR) if $h'(t) < 0$ for $t \in (0, t_0)$, $h'(t_0) = 0$ and $h'(t) > 0$ for $t > t_0$, for a comprehensive account of bathtub failure rate distributions see Rajarshi and Rajarshi (1988);
4. Upside-down bathtub shaped failure rate (UBFR) if $h'(t) > 0$ for $t \in (0, t_0)$, $h'(t_0) = 0$ and $h'(t) < 0$ for $t > t_0$;
5. Modified bathtub shaped failure rate (MBFR) if $h(t)$ is first increasing and then bathtub shaped;
6. Roller-coaster shaped failure rate if there exist *n* consecutive change points $0 < t_1 < t_2 < ... < t_n < \infty$ such that in each interval $[t_{j-1}, t_j]$, $1 \leq j \leq n+1$, where $t_0 = 0, t_{n+1} = \infty$, $h(t)$ is strictly monotone and it has opposite monotonicity in any two adjacent such intervals.

**Definition 2.2** (Lariviere and Porteus (2001))

Let *X* be a nonnegative random variable, then the generalized failure rate of *X* is *g(x) = x h(x)*.

If *g(x)* is increasing weakly in *x* for *F(x) < 1*, then the random variable *X* is said to have an increasing generalized failure rate (IGFR) or *F(x)* is an IGFR distribution. Similarly, decreasing the generalized failure rate (DGFR) can be defined. IGFR distributions are mainly considered in pricing and supply chain management.

**Theorem 2.1** (Lariviere (2006))

The following statements are equivalent:

1. *X* is IGFR.
2. *log(X)* is IFR.
3. $X \leq_{hr} \lambda X$ for $\lambda \geq 1$,
4. $f(x, \theta) = \overline{F}(x/\theta)$ is TP$_2$

Where $X_1 \leq_{hr} X_2$ denotes $X_1$ is smaller than $X_2$ in the hazard rate order (Rose *et al.*(1996)).

**Table 2: Preservation of Life Distributions under Reliability Operations**

| Class | Reliability operation | | |
|---|---|---|---|
| | **Coherent Systems** | **Convolution** | **Mixture** |
| IFR | Not Preserve | Preserve | Not Preserve |
| IFRA | Preserve | Preserve | Not Preserve |
| NBU | Preserve | Preserve | Not Preserve |
| NBUE | Not Preserve | Preserve | Not Preserve |
| DMRL | Not Preserve | Not Preserve | Not Preserve |
| HNBUE | Not Preserve | Preserve | Not Preserve |
| NBU-$t_0$ | Preserve | Not Preserve | Not Preserve |
| DFR | Not Preserve | Not Preserve | Preserve |
| DFRA | Not Preserve | Not Preserve | Preserve |
| NWU | Not Preserve | Not Preserve | Not Preserve |
| NWUE | Not Preserve | Not Preserve | Not Preserve |
| IMRL | Not Preserve | Not Preserve | Preserve |
| HNWUE | Not Preserve | Not Preserve | Preserve |
| NWU-$t_0$ | Not Preserve | Not Preserve | Not Preserve |
| BFR | Not Preserve | Not Preserve | Not Preserve |

*1.1 Properties of Ageing classes*

Some properties of ageing distributions can be found in Barlow and Campo (1975), Lee and Thompson (1976), Bergman (1977), Langberg *et al.*(1980), Klefsjo (1982), Aarset (1987) and Xie (1989). A few of them are listed below.

Suppose $\phi(p)$ be the scaled TTT transform of a continuous lifetime distribution *F*, then

1. *F* is IFR (DFR) if and only if $\phi(p)$ is concave (convex) in $p \in [0,1]$.

2. $F$ is IFRA (DFRA) if and only if $\phi(p)/p$ is decreasing (increasing) in $p \in [0,1]$.

3. $F$ is NBUE (NWUE) if and only if $\phi(p) \geq p$ $(\phi(p) \leq p)$ for $p \in [0,1]$.

4. $F$ is DMRL (IMRL) if and only if $(1-\phi(p))/(1-p)$ is decreasing (increasing) in $p \in [0,1]$.

5. $F$ is HNBUE (HNWUE) if and only if $\phi(p) \leq 1 - \exp\{-F^{-1}(p)/\mu\}$ $(\phi(p) \geq 1 - \exp\{-F^{-1}(p)/\mu\})$ for $p \in [0,1]$.

6. $F \in$ BFR (UFR) if $\phi$ has only reflection point $u_0$ such that $0 < u_0 < 1$ and it is (convex) on $[0, u_0]$ and concave (convex) on $[u_0, 1]$.

Table 2 illustrate that the life distributions are preserved under various reliability operations.

## 3. IFR and DFR Distributions

Distributions with monotone failure rates have been considered a significant factor because of their practical interest in reliability. Such distributions comprise a massive class. This section deals with mainly the characteristics of IFR and DFR class distributions.

*3.1 Increasing Failure Rate*

In the class of failure rate distributions, increasing failure rate (IFR) distributions have their relevance. (see, for example, Marshall and Proschan (1965), Kvam (2002) and Koutras (2011)). The Weibull and Gamma are two primer distributions that are effectively using in many industries for optimizing system reliability, especially for redundancy allocation problems (RAP) (see, Fyffe *et al.*(1968)).

If $F$ is IFR, it is equivalent to say that $F$ is Polya frequency function of order 2 (PF2). Schoenberg (1951) introduced PF2 and defined by, a PF2 function is a non-negative measurable function $g(x)$ defined for all real $x$ such that

$$\begin{vmatrix} g(x_1 - y_1) & g(x_1 - y_2) \\ g(x_2 - y_1) & g(x_2 - y_2) \end{vmatrix} \geq 0,$$

whenever $x_1 < x_2$ and $y_1 < y_2$; in addition, $g(x) \neq 0$ for at least two distinct values of $x$. If the density function $f$ has PF2 then it would be easy to verify that the distribution has IFR, but the converse is needed not to be true. PF2 has many useful properties including unimodality, closure under convolution, variation diminishing and certain moment properties (see Schoenberg (1951) and Karlin (1961)).

**Theorem 3.1.1** (Barlow, Marshall, Proschan (1963))

If F and G are IFR, then their convolution H, given by $H(t) = \int_{-\infty}^{\infty} F(t-x) dG(x)$ is also IFR.

Bakouch *et al.* (2014) proposed a model (Binomial-exponential 2) using zero-truncated binomial random variable having IFR,

$$h(x) = \lambda \left(1 - \frac{\theta}{2 - \theta + \lambda \theta x}\right).$$

It is increasing in *x* for $0 \leq \theta \leq 1$ and $\lambda > 0$.

Guilani *et al.* (2016) investigated RAP of a series-parallel system with components having IFR based on Weibull distribution.

Bugatekin (2017) developed a two-parameter mixed distribution named RL (Rayleigh-Logarithmic) distribution with IFR function,

$$h(x) = -\frac{xe^{\frac{-x^2}{2\sigma^2}} p \left(1 - e^{\frac{-x^2}{2\sigma^2}} p\right)^{-1}}{\sigma^2 \ln\left(1 - e^{\frac{-x^2}{2\sigma^2}} p\right)}, \sigma > 0, y \in [0, \infty).$$

Elbarmi *et al.* (2020) provide IFRA model with an estimator of *F* that is uniformly strongly consistent and they derive convergence of the estimator at the point where *F* is IFRA using the arg max theorem.

*3.2 Decreasing Failure Rate*

It is reasonable to assume that the failure rate of the life distribution is increasing (IFR). The wear-outphase easily interprets this fact, implying that the age-oldunits have a higher chance of failure. However, it appears to be more difficult to explain why a unit's life expectancy is increasing while its failure rate is decreasing. Decreasing failure rate (DFR) distributions may appear in different ways including, when the strength of some metals increases with usage or when an instrument exhibit an 'infant mortality' rate during their lifetime. Normally DFR represents the apparatus having betterment with age so that age-old units have less chance of failure.

**Theorem 3.2.1** (Barlow, Marshall, Proschan (1963))

- $F$ is DFR if and only if the support of $F$ is $[0, \infty)$, and $1 - F(x+y)$ is TP2 for $x + y \geqq 0$.
- DFR property is not preserved under the reliability conditions of convolution and coherent systems. However, mixtures of DFR distributions are DFR.
- If $F(t, \phi)$ is a DFR distribution in $t$ for each $\phi$ in $\Phi$, then $G(t) = \int_\Phi F(t, \phi) d\mu(\phi)$ is DFR where $\mu$ is a probability measure in $\Phi$.

Adamidis and Loukas (1998) presented a DFR distribution named Exponential–Geometric distribution obtained as a mixture of the exponential and geometric distribution. Its failure rate $h(x) = \beta(1 - pe^{-\beta x})^{-1}, p \in (0,1), \beta > 0, x > 0$ is decreasing in $x$.

Finkelstein and Esaulova (2001) observed a mixture of several continuous IFR distributions and analyzed the limiting behaviour of the mixture failure rate function. They found the mixture can be DFR and for that the limiting behaviour of conditional mean and the conditional variance of the mixing parameters are essential.

Chahkandi and Ganjali (2009) introduced a DFR distribution which is a mixture of power-series distribution and exponential distribution. Generally, mixtures of exponential distributions are DFR.

*3.3 Basic Properties of IFR and DFR distributions*

Pham (2003) provides some properties of IFR and DFR distributions.

- Suppose $X_1$ and $X_2$ are IFR, then $X_1 + X_2$ are also IFR; but DFR property may not be preserved. In the case of mixtures, the opposite happens.
- IFR distributions are preserved under coherent systems.
- Order statistics of an IFR distribution are again IFR; but not necessarily true for DFR distributions.
- Spacings of DFR distributions gives DFR; but not true in the case of IFR distributions
- The *pdf* of IFR distributions may not be unimodal but a decreasing function for DFR distributions.

## 4. Bathtub and Upside-down Bathtub Shaped Failure Rate

In most real-life situations failure rates may not exhibit monotonic nature rather they start with a peak value, decrease, may then persist constant and in the end increases rapidly (see, Rajarshi and Rajarshi (1988)). These indicate the different aspects of a component or system, *i.e.*, early life (infant mortality phase), useful life and wear-out phase. They have also known, respectively as burn-in, random and wear out phases in the reliability context. This occurs mainly due to diverse amplitudes of the life of a device when they are concerned by a collection of deformities. In modelling real-life data sets, bathtub failure rates are beneficial because the lifetimes of mechanical products, electromechanical and electronics are frequently modelled with this aspect.

*4.1 Definitions and results for BFR distributions*

**Definition 4.1.1** (Glaser (1980))

Suppose $F$ is a lifetime distribution that is differentiable and having a density $f$ and failure rate $\lambda$. Then $F$ is said to be a BFR distribution, if there exists a $t_o$, $0 < t_o < \infty$ such that $\lambda(t)$ is decreasing for $t \leq t_o$ and increasing for $t_o \leq t$.

**Definition 4.1.2** (Deshpande and Suresh (1990))

A life distribution $F$ is said to be a BFR if there is a point $t_0 > 0$ for which $-\log \overline{F}(t)$ is concave in $[0, t_o]$ and convex in $[t_0, \infty)$.

**Definition 4.1.3** (Deshpande and Suresh (1990))

A life distribution is known to be an "increasing initially, then decreasing mean residual life" (IDMRL) distribution if there is a point $t_o > 0$ such that

- $\mu(s) \leq \mu(t), 0 \leq s < t < t_0$ and

- $\mu(u) \geq \mu(v), t_0 \leq u < v < \infty$.

IDMRL distributions are also known as upside-down bathtub shaped mean residual life (UBTMRL) distribution (see Rajarshi and Rajarshi (1988)).

**Definition 4.1.4** (Mi (1995))

A distribution function $F$ is said to be a BFR (UBFR) type if there exist two points $t_1$ and $t_2$ for which $0 \leq t_1 \leq t_2 < \infty$ such that the failure rate function $h(t)$ is,

1. strictly decreasing (increasing) in $0 \leq t \leq t_1$;
2. constant in $t_1 \leq t \leq t_2$; and
3. strictly increasing (decreasing) for $t_2 \leq t$.

where $t_1$ and $t_2$ are called change points of $h(t)$. BFR turn out to strictly increasing (IFR) when $t_1 = t_2 = 0$. If $t_1 = t_2 \to \infty$, $h(t)$ becomes strictly decreasing (DFR). In common, when $t_1 = t_2$ the interval for which $h(t)$ constant degenerates to a single point. This definition treats monotonic failure rates IFR and DFR as a special case of BFR.

**Definition 4.1.5** (Mitra and Basu(1995))

An *absolutely continuous distribution* with support $[0, \infty)$ is called a BFR distribution if its failure rate function, $h(t)$ is non-increasing in $[0, t_0)$ and non-decreasing in $[t_0, \infty)$ for a $t_0 \geq 0$.

In definition, a maximum of two change points is permitted. In line with Mi (1995) the points in the interval $(t_1, t_2)$ are not change points but they appeared as change points in Mitra and Basu (1995).

A generalization of definition 1 is given by,

**Definition 4.1.6** (Haupt and Scabe (1997))

A life distribution $F$ is said to be BFR, if there is a point $t_o$, $0 < t_o < \infty$, and $\overline{F}(x|t) = \overline{F}(t+x)/\overline{F}(t)$, for which

1. For $0 \leq t \leq t_0$, and $0 \leq x \leq t_0 - t$, $\overline{F}(x|t)$ is increasing w.r.t. $t$; and

2. For $t_0 \leq t \leq t_\infty$ and $x \geq 0$, $\overline{F}(x|t)$ is decreasing w.r.t. $t$.

**Theorem 4.1.1** (Olcay (1995))

If $h(t)$ is the failure rate of a BFR distribution, then

1) $\mu(t)$ is decreasing if $h(0) \leq 1/\mu$.

2) $\mu(t)$ is of the type UBFR if $h(0) > 1/\mu$.

This theorem examines the behaviour of MRL property for BFR type distributions. Similarly, the succeeding theorem validates the behaviour of UBFR type distributions.

**Theorem 4.1.2** (Olcay (1995))

Suppose $h(t)$ is of type UBFR; then

1) $\mu(t)$ is increasing if $h(0) \geq 1/\mu$.
2) $\mu(t)$ is of type BFR if $h(0) < 1/\mu$.

*4.2 Properties of BFR distributions*

Mitra and Basu (1996) provided various properties of BFR distributions. Preservations of lifetime distributions under the reliability operations like convolution, mixtures and coherent systems are viewed.

- Suppose $F$ is BFR and $G$ is exponential with mean $\{h(t_0)\}^{-1}$, then $\overline{F}(t) \leq \overline{G}(t)$, where $t_o$ is the change point for $h(t)$ is minimum.
- $E(X^k) \leq \dfrac{\Gamma(k+1)}{r(t_0)^k}, k > 0$.
- If a BFR distribution has a $k^{th}$ raw moment, $E(X^k) = \dfrac{\Gamma(k+1)}{\{r(t_0)\}^k}$ then it is necessarily exponential.
- If the BFR distribution has a finite turning point $t_0$, then $\{\overline{F}(x)/\overline{F}(t_0)\}^{1/(x-t_0)}$ is decreasing in $x, x \in (t_0, \infty)$.
- The moment sequence uniquely determines the distribution if it is BFR with a finite turning point.
- Convolution of BFR distributions is not necessarily BFR.
- A mixture of BFR distributions may not be BFR.
- A parallel system consisting of two independent components from BFR distributions may not be BFR.

*4.3 Bathtub construction techniques*

Various techniques are used to construct BFR type distributions. A few of them are given.

- **Glaser's technique**

    Glaser (1980) adopted a function $\eta(t)$ that accomplishes the succeeding norms:

a) $\eta(t) = -f'(t)/f(t)$ and $f(t)$ is a density function;

b) there exist a $t_o > 0$ such that $\eta'(t) < 0$ for all $t \in (0, t_o)$, $\eta'(t_0) = 0$ and $\eta'(t) > 0$ for all $t > 0$;

c) there exist a $y_o > 0$ such that $\int_{y_0}^{\infty} \left[ f(y)/f(y_0) \right] \eta(y_0) dy - 1 = 0$.

- **Convex function**

    From the description of BFR distribution, we can construct a BFR by selecting a positive convex function $h(t)$ that satisfying $\int_{0}^{\infty} h(t) dt = \infty$ (Rajarshi and Rajarshi (1988)).

- **Function of random variables**

    This technique due to Griffith (1982), are for differentiable failure rate functions.

- **Reliability and stochastic Mechanisms**

    Series system (competing risk model). A model gained by Murthy *et al.* (1973)

- **Mixtures**

    Mixtures of BFR distributions frequently produces BFR distributions.

- **Truncating DFR models**

    BFR models can be constructed by truncating DFR type distributions.

*4.4 T*otal time on test *transform*

The notion of total time on test (TTT-transform) transform is very popular given its relevancy in diverse areas including reliability analysis, stochastic modelling, ordering distributions, econometrics, *etc*. TTT transform is found to be very useful for the study of ageing properties also used to construct BFR distributions.

The TTT-transform of a life distribution F is defined as

$$H_F^{-1}(t) = \int_{0}^{F^{-1}(t)} \overline{F}(u) du, 0 \leq t \leq 1.$$

If $\mu$ represents the mean lifetime, then $H_F^{-1}(1) = \int_0^{F^{-1}(1)} \overline{F}(x)\,dx$, where $\mu$ is finite.

The term $\phi(p)$, where $\phi(p) = H_F^{-1}(p)/H_F^{-1}(1) = H_F^{-1}(p)/\mu$ is known as scaled TTT transform. The equilibrium distribution of $F$ when F is non-arithmetic is the same as the distribution of $\phi(p)$ for $p \in [0,1]$.

A generalization of scaled TTT transform introduced Jewell (1977). Bergman (1979) proposed a method for identifying exponentiality versus BFR based on total time on the test. For the detailed study of major properties and applications of TTT transform see Li and Zou (2004), Lai and Xie (2006), Li and Shaked (2007), Nair *et al*. (2008) and Esfahani (2021). Increasing convex (concave) TTT transform is examined by Deepthi and Chacko (2021) to identify the failure rates of functions of random variables.

Examples:

- Chacko (2016) introduced a BFR model named an *x*-Exponential model similar to Generalized Lindley with failure rate,

$$h(x) = \frac{\alpha\left(1 - (1+\lambda x^2)e^{-\lambda x}\right)^{\alpha-1}\left[\lambda x^2 - 2x + 1\right]\lambda e^{-\lambda x}}{1 - \left(1 - (1+\lambda x^2)e^{-\lambda x}\right)^{\alpha}}, x > 0, \alpha > 0, \lambda > 0.$$

- Arshad *et. al* (2021) developed a flexible lifetime model with BFR characteristics, its failure rate is given by,

$$h(x) = \frac{\beta(\alpha+1)}{\alpha g\left(1 - \dfrac{x}{g}\right)\left(1 + \dfrac{x}{\alpha g}\right)}, \text{ where } 0 < x \leq g, \alpha > 0.$$

## 5. Roller-Coaster Shaped Failure Rate

Often some distributions arise with more than one turning point for the failure rates. Roller-coaster failure rates appear in such situations (see Wong (1988), Wong and Lindstrom (1988) and Wong (1991)). The monotonicity of failure rates with one turning point can be easily determined using Glaser's technique. To study the monotonicity of failure rates with more than one turning point, Gupta and Warren (2001) introduced generalized Glaser's technique and examined several examples. Moreover, they found out the relation between turning points of failure rates and Glaser's eta function and put forward the number of turning points of failure rates that do not surpass the number of turning points of Glaser's eta function.

Gupta and Viles (2011) have taken on the monotonicity of Glaser's eta function, failure rate and mean residual life function and the results are used to investigate the results in extended generalized inverse Gaussian distribution.

Table 3 illustrates some of the recent lifetime distributions and their failure rate classes.

**Table 3 Distributions and their Failure rates**

| Distribution | Failure rates | Reference |
| --- | --- | --- |
| Generalized Exponential Geometric | IFR,DFR,UBFR | Silva *et al.* (2010) |
| Exponential-geometric range | IFR | Shahsanaei *et al.* (2012) |
| Exponential Poisson-Lindley | DFR | Barreto-Souza and Bakouch (2013) |
| Nadarajah–Haghighi | Constant,IFR,DFR,BFR,UFR | Lemonte (2013) |
| Kumaraswamy generalized Rayleigh | IFR,DFR,BFR | Gomes *et al.* (2014) |
| Exponentiated Exponential-Geometric | IFR,DFR,UBFR | Louzada *et al.* (2014) |
| Marshall-Olkin Generalized Exponential | IFR, DFR, BFR,UBFR | Ristić and Kundu (2015) |
| Extended Inverse Lindley | UBFR | Alkarni(2015) |
| Modified weibull geometric | IFR,DFR,BFR,UBFR | Wang and Elbatal (2015) |

| Generalized Bilal | IFR,DFR,UBFR | Abd-Elrahman(2017) |
| --- | --- | --- |
| Alpha Power-transformed Weibull | Constant, IFR, DFR, BFR, UBFR, MFR | Dey (2017) |
| Generalized inverse xgamma | IFR,DFR,UBFR | Tripathi *et al.* (2018) |
| Generalized Weibull Uniform | IFR,DFR,BFR | Khaleel *et al.* (2019) |
| Generalized X-Exponential | IFR,DFR,BFR | Chacko and Deepthi(2019) |
| Marshall-Olkin logistic-exponential | IFR,DFR,BFR,UBFR | Mansoor *et al*. (2019) |
| Burr–Hatke Exponential | DFR | Yadav *et al*. (2019) |
| Odd Lindley-Inverse exponential | DFR | Ieren and Abdullahi (2020) |
| Generalized Lindley | IFR,DFR,BFR | Maurya *et al.* (2020) |
| Inverted Power Rama | DFR,UBFR | Onyekwere *et al.* (2020) |
| Generalized Log-Weibull | IFR,DFR,BFR,UBFR, S-shape | Kumar and Nair (2021) |

## 6. Generalizations of Lindley Distribution

Recently Lindley distribution was introduced by Lindley (1958) recognized for its applications in studying lifetime data. Various generalizations of Lindley distribution are available in the statistical literature. Here some of the generalizations and their failure rate shapes have been discussed.

- Nadarajah *et al.* (2011) introduced two parameters generalized Lindley distribution with failure rate,

$$h(x) = \frac{\alpha\lambda^2}{1+\lambda}(1+x)\left[1-\frac{1+\lambda+\lambda x}{1+\lambda}\exp(-\lambda x)\right]^{\alpha-1}\exp(-\lambda x)\{1-V^\alpha(x)\}^{-1},$$

where $V(x) = 1 - \frac{1+\lambda+\lambda x}{1+\lambda}\exp(-\lambda x)$, $\alpha > 0$ and $x > 0$.

This distribution has DFR and BFR for $\alpha < 1$ and IFR for $\alpha \geq 1$.

- Bakouch *et al.* (2012) introduced an extended version of Lindley distribution (Extended Lindley distribution) by considering particular exponentiation of Lindley distribution proposed by Lindley (1958). Its failure rate is,

$$h(x) = \frac{\beta(1+\lambda+\lambda x)\lambda^\beta x^{\beta-1} - \lambda\alpha}{1+\lambda+\lambda x}, x > 0, \alpha, \beta, \lambda > 0.$$

The distribution has IFR for $\alpha > k$ and DFR for $\alpha < k$, where $k = -\beta(\beta-1)(\lambda x)^{\beta-2}(1+\lambda+\lambda x)^2$. For some other combinations of parameter values, it behaves as BFR and UBFR distributions.

- Ibrahim *et al.* (2013) proposed a new generalized Lindley distribution, obtained from a mixture of Gamma($\beta,\theta$) and Gamma($\beta,\theta$) with mixing constant $\frac{\theta}{1+\theta}$. The failure rate function is,

$$h(x) = \frac{\frac{1}{1+\theta}\left[\frac{\theta^{\alpha+1}x^{\alpha-1}}{\Gamma(\alpha)} + \frac{\theta^\beta x^{\beta-1}}{\Gamma(\beta)}\right]e^{-\theta x}}{1 - \frac{1}{1+\theta}\left[\frac{\theta\gamma(\alpha,\theta x)}{\Gamma(\alpha)} + \frac{\gamma(\beta,\theta x)}{\Gamma(\beta)}\right]}; \alpha, \theta > 0, x > 0.$$

It shows IFR, DFR, BFR, UBFR and modified bathtub shaped failure rates.

- Beta-generalized Lindley (4-parameter class of generalized Lindley) distribution proposed by Oluyede and Yang (2015) with density function,

$$h(x) = \frac{\alpha\lambda^2}{B(a,b)(1+\lambda)}(1+x)\exp(-\lambda x)\left[1 - \frac{1+\lambda+\lambda x}{1+\lambda}\exp(-\lambda x)\right]^{a\alpha-1}$$

$$\times \left\{ 1 - \left[ 1 - \frac{1+\lambda+\lambda x}{1+\lambda} \exp(-\lambda x) \right]^{\alpha} \right\}^{b-1} \left[ 1 - \frac{1}{B(a,b)} \int_{0}^{G(x;\lambda,\alpha)} t^{a-1}(1-t)^{b-1} \right]^{-1} ; x > 0,$$

$\alpha, \lambda, a, b > 0$. Where $G(x;\alpha,\lambda) = \left[ 1 - \frac{1+\lambda+\lambda x}{1+\lambda} \exp(-\lambda x) \right]^{\alpha}$.

The distribution has IFR, DFR and BFR shape properties for different choice of parameters.

- Al-Babtain *et al.* (2015) introduced five parameters Lindley distribution, a generalization of the basic Lindley distribution. It has the failure rate,

$$h(x) = \frac{\theta^2 \left[ \frac{k(\theta x)^{\alpha-1}}{\Gamma(\alpha)} + \frac{\eta(\theta x)^{\beta-1}}{\theta \Gamma(\beta)} \right] e^{-\theta x}}{\eta + \theta k - \left[ \theta k \gamma_\alpha(\theta x) + \eta \gamma_\beta(\theta x) \right]}, x > 0, \theta, \alpha, \beta > 0, k, \eta \geq 0$$

where $\gamma_a(b) = \frac{\gamma(a,b)}{\Gamma(a)} = \frac{1}{\Gamma(a)} \int_{0}^{b} t^{a-1} e^{-t} dt$. Depending on the parameters, its failure rate might have constant, increasing, decreasing or bathtub shaped.

- Oluyede *et al.* (2015) proposed log generalized Lindley-Weibull distribution which is a generalization of the Lindley distribution via the Weibull model. The failure rate of the given distribution is,

$$h(x) = \frac{\frac{c\theta^{\alpha+1}}{\gamma(\beta+\theta)\Gamma(\alpha+1)} \left(\frac{x}{\gamma}\right)^{c\alpha-1} \left\{ \alpha + \beta \left(\frac{x}{\gamma}\right)^c \right\} e^{-\theta\left(\frac{x}{\gamma}\right)^c}}{1 - \frac{1}{(\beta+\theta)\Gamma(\alpha)} \left\{ \theta \left[ \Gamma(\alpha) - \Gamma(\alpha,u) \right] + \frac{\beta}{\alpha} \left[ \Gamma(\alpha) - \Gamma(\alpha+1,u) \right] \right\}};$$

$\theta, \alpha, c, \gamma, \beta > 0$ and $x > 0$. The distribution might have IFR, DFR or BFR depending on the parameters.

- Bhati *et al.* (2016) introduced a new 3-parameter extension of generalized Lindley distribution with failure rate,

$$h(x) = \frac{\alpha \theta^2 \lambda (1+x\lambda)^{-1+2\alpha}}{1+\theta(1+x\lambda)^{\alpha}}; \alpha, \lambda, \theta > 0, x > 0.$$

The failure rate function is decreasing, upside down and increasing according to as

$$\left(0<\alpha<\frac{1}{2}\right)\cup\left(\frac{1}{2}<\alpha<1\cap\theta>\frac{1-2\alpha}{\alpha(\alpha-1)}\right), \quad \left(\frac{1}{2}<\alpha<1\cap\theta>\frac{1-2\alpha}{\alpha(\alpha-1)}\right) \quad \text{and}$$

$(\alpha>1)\cap(\theta>0)$ respectively.

- Based on modifications of the quasi Lindley distribution Elgarhy (2016) introduced transmuted generalized Lindley distribution. It has the density function,

$$h(x)=\frac{\frac{a\theta^2}{\theta+1}(1+x)e^{-\theta x}\left[1-e^{-\theta x}\left(1+\frac{\theta x}{\theta+1}\right)\right]^{a-1}\left\{(1+\lambda)-2\lambda\left[1-e^{-\theta x}\left(1+\frac{\theta x}{\theta+1}\right)\right]^a\right\}}{1-\left[1-e^{-\theta x}\left(1+\frac{\theta x}{\theta+1}\right)\right]^a\left\{1+\lambda-\lambda\left[1-e^{-\theta x}\left(1+\frac{\theta x}{\theta+1}\right)\right]^a\right\}}; \text{ for}$$

$x>0, a>0, \theta>-1$. The distribution has IFR for all values of parameters.

- Ramos and Louzada (2016) proposed three-parameter generalized weighted Lindley distribution with failure rate,

$$h(x)=\frac{\alpha\lambda^{\alpha\phi}x^{\alpha\phi-1}\left(\lambda+(\lambda x)^\alpha\right)e^{-(\lambda x)^\alpha}}{\Gamma\left[\phi,(\lambda x)^\alpha\right](\lambda+\phi)+(\lambda x)^{\alpha\phi}e^{-(\lambda x)^\alpha}}, x>0, \alpha,\lambda,\phi>0.$$

The proposed distribution has IFR, DFR, BFR, UBFR or decreasing-increasing-decreasing shape properties.

- Two-parameter generalized inverse Lindley distribution introduced Sharma *et al.* (2016). The failure rate function is unimodal and is given by,

$$h(x)=\frac{\alpha\theta^2\left(1+x^\alpha\right)}{x^{\alpha+1}\left[(1+\theta)x^\alpha\left(e^{\frac{\theta}{x^\alpha}}-1\right)-\theta\right]}, x>0, \alpha,\theta>0.$$

- Benkhelifa (2017) defined a new 3-parameter model called the Marshall-Olkin extended generalized Lindley distribution with failure rate,

$$h(x) = \frac{\alpha\theta^2 e^{-\theta x}(1+x)\left(1-\frac{1+\theta+\theta x}{1+\theta}e^{-\theta x}\right)^{\alpha-1}}{(1+\theta)\left[\beta+\bar{\beta}\left(1-\frac{1+\theta+\theta x}{1+\theta}e^{-\theta x}\right)^{\alpha}\right]\left[1-\left(1-\frac{1+\theta+\theta x}{1+\theta}e^{-\theta x}\right)^{\alpha}\right]},$$

for $\alpha, \beta, \theta > 0$ and $x > 0$. The failure rate function can be increasing, decreasing, upside-down bathtub (unimodal), bathtub-shaped or increasing-decreasing-increasing depending on the parameter values.

- Ekhosuehi and Opone (2018) proposed a three-parameter generalized Lindley distribution with failure rate,

$$h(x) = \frac{\alpha\lambda^2(\beta+x^\alpha)x^{\alpha-1}}{(1+\lambda\beta+\lambda x^\alpha)}, x > 0, \theta, \lambda, \beta > 0.$$

For $\alpha < 1$, it has DFR and IFR for $\alpha \geq 1$.

- Abouammoh and Kayid (2020) presented a new method for generalizing Lindley distribution, by increasing the number of mixed models. The introduced model generalized Lindley of order $m$ has the failure rate,

$$h(x) = \frac{\theta^m \sum_{i=1}^{m} \frac{x^{m-i}}{\Gamma(m-i+1)}}{\sum_{j=0}^{m-1}\sum_{i=1}^{m-j} \theta^{m-i} \frac{x^j}{j!}}, \theta > 0, x \geq 0.$$

The distribution has IFR for all values $m$ and $\theta$.

- Algarni (2021) proposed an extension of generalized Lindley distribution using Marshall-Olkin method called M-O Extended Generalized Lindley distribution. Its failure rate is,

$$h(x) = \frac{\lambda^2 \gamma}{(\lambda+1)} \frac{(x+1)K^{\gamma-1}e^{-\lambda t}}{(K^\gamma - 1)(\delta(K^\gamma - 1) - K^\gamma)},$$

where $K = 1 - \frac{(\lambda + \lambda t + 1)e^{-\lambda t}}{\lambda + 1}$, $\lambda, \gamma, \delta > 0$ and $x > 0$. The distribution has IFR, DFR and BFR shape properties.

Even though the Lindley distribution belongs to the class of increasing failure rates, many of its generalizations raise it to different classes. Most of the generalizations bring flexibility to the distribution.

## 6. DUS Transformation

In statistical literature, many different kinds of transformations have been introduced to take the distribution to another level. Generally, transformation allows flexibility to the baseline distribution. Kumar *et al.* (2015) introduced a transformation method called DUS transformation, which gives parsimonious distribution in terms of computation and interpretation. Since it doesn't add a new parameter, it is clearly a transformation and not a generalization. According to Kumar *et al.* (2015), if $f(x)$ and $F(x)$ be the *pdf* and *cdf* of some baseline equation, then the *pdf* and *cdf* of DUS transformed distribution are $g(x) = \frac{1}{e-1} f(x) e^{F(x)}$ and $G(x) = \frac{1}{e-1}\left[e^{F(x)} - 1\right]$, respectively. They considered exponential distribution as the baseline distribution and the resultant gives an IFR distribution even the exponential distribution has only a constant failure rate. The failure rate of the DUS-exponential distribution is $h(x) = \theta e^{-\theta x}\left[e^{e^{-\theta x}} - 1\right]^{-1}, x > 0, \theta > 0$.

Maurya *et al.* (2017) considered Lindley distribution as the baseline distribution and named it exponential transformed Lindley distribution. It has a failure rate $h(x) = \frac{\theta^2 (1+x) e^{-\theta x}}{(\theta+1)\left(\exp\left[e^{\theta x} \frac{1+\theta+\theta x}{1+\theta}\right] - 1\right)}, x > 0, \theta > 0$. The distribution has IFR for all the values of $\theta$.

DUS transformation of Lomax distribution (DUS-Lomax distribution) was introduced by Deepthi and Chacko (2020). Its failure rate is given by, $h(x) = \alpha\beta(1+\beta x)^{-(\alpha+1)}\left[e^{(1+\beta x)^{-\alpha}} - 1\right]^{-1}, x > 0, \alpha, \beta > 0.$

DUS-Lomax has DFR and UBFR shape properties even though the Lomax distribution belongs to the family of DFR.

A generalized lifetime model has been proposed by Kavya and Manoharan (2020) using generalized DUS transformation with Weibull distribution as the baseline distribution (GDUS-Weibull distribution). The Weibull distribution has only monotone failure rates but GDUS-

Weibull distribution has IFR, DFR and UBFR shape properties, and the failure rate is given by,

$$h(x) = \alpha k x^{k-1} e^{-\left(\frac{x}{\lambda}\right)^k} \left(1 - e^{-\left(\frac{x}{\lambda}\right)^k}\right)^{\alpha-1} e^{\left(1 - e^{-\left(\frac{x}{\lambda}\right)^k}\right)^\alpha} \lambda^{-k} \left(e - e^{\left(1 - e^{-\left(\frac{x}{\lambda}\right)^k}\right)^\alpha}\right)^{-1}, x > 0, \alpha, \lambda, k > 0.$$

Karakaya *et al.* (2021) considered Kumaraswamy distribution as the baseline distribution with failure rate,

$$h(x) = \frac{\alpha \beta x^{\alpha-1} (1-x^\alpha)^{\beta-1} e^{\left(1-(1-x^\alpha)^\beta\right)}}{e - e^{\left(1-(1-x^\alpha)^\beta\right)}}, 0 < x < 1, \alpha, \beta > 0.$$

The proposed distribution has IFR and BFR.

DUS transformation of inverse Weibull distribution has been proposed by Gauthami and Chacko (2021). The distribution has DFR and UBFR and its failure rate is

$$h(x) = \frac{1}{e - e^{e^{-\left(\frac{x}{\beta}\right)^{-\alpha}}}} \left(\frac{\alpha}{\beta}\right) \left(\frac{x}{\beta}\right)^{-(\alpha+1)} e^{\left\{-\left(\frac{x}{\beta}\right)^{-\alpha} + e^{-\left(\frac{x}{\beta}\right)^{-\alpha}}\right\}}, x > 0, \alpha, \beta > 0.$$

In studying the above resolutions, DUS transformation plays a crucial role in adding flexibility in the baseline distributions and they have been considering more suitable for modelling lifetime data.

*6.1 DUS Exponential-Weibull distribution*

In this section, a new failure rate distribution using DUS transformation of the mixture of Exponential and Weibull distribution introduced by Anakha and Chacko (2020), is proposed. The *pdf* of the Expoential-Weibull distribution is,

$$f(x) = \frac{\lambda^2 e^{-\lambda x} + \alpha \lambda^\alpha x^{\alpha-1} e^{-(\lambda x)^\alpha}}{1+\lambda}, \alpha, \lambda > 0, x > 0.$$

Then *pdf* of DUS Exponential-Weibull distribution, denoted by DUS-EW, is obtained as,

$$g(x) = \frac{\lambda^2 e^{-\lambda x} + \alpha \lambda^\alpha x^{\alpha-1} e^{-(\lambda x)^\alpha}}{(e-1)(1+\lambda)} e^{1-V(x)}, \text{ where } V(x) = \frac{\lambda e^{-\lambda x} + e^{-(\lambda x)^\alpha}}{1+\lambda}.$$

The associated failure rate is,

$$h(x) = \frac{\lambda^2 e^{-\lambda x} + \alpha \lambda^\alpha x^{\alpha-1} e^{-(\lambda x)^\alpha}}{\left(e - e^{1-V(x)}\right)(1+\lambda)} e^{1-V(x)}, \alpha, \lambda > 0, x > 0.$$

The DUS-EW distribution has increasing, decreasing and and upside-down bathtub shaped failure rate properties; see Figure 6.1.

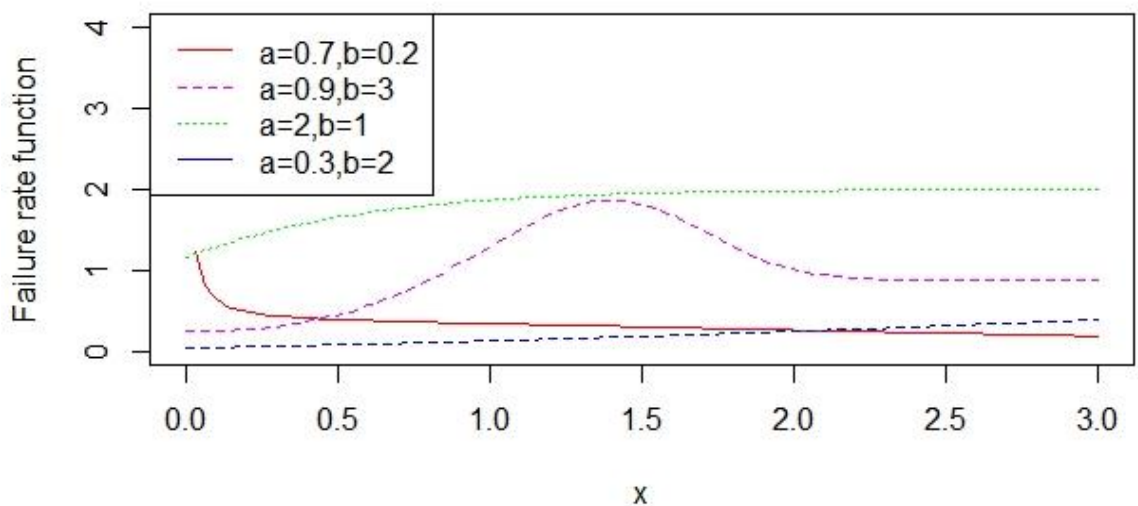

**Figure 6.1 Failure rate function of DUS Exponential-Weibull distribution**

7. **Conclusion**

This paper examined several characteristics of ageing properties and different shapes of failure rates of some recent lifetime distributions. Distributions with different failure rate shape properties are more favourable in structuring probability models. Various generalizations of Lindley distribution are discussed and the distinguishable roles of DUS transformation are discussed. A new failure rate distribution using DUS transformation of the mixture of Exponential and Weibull distribution is proposed.